\documentclass{ifacconf}
\usepackage{natbib}
\usepackage{graphicx}
\usepackage{amssymb}
\usepackage[cmex10]{amsmath}
\usepackage{array}
\usepackage{mdwmath}
\usepackage{mdwtab}
\usepackage{balance}
\usepackage[usenames, dvipsnames]{color}
\usepackage{epstopdf}
\usepackage{tikz}

\usetikzlibrary{shapes,arrows}

\newtheorem{theorem}{Theorem}
\newtheorem{definition}{Definition}
\newtheorem{lemma}{Lemma}

\renewcommand\Re{\operatorname{Re}}

\newcommand{\ie}{\textit{i}.\textit{e}.}

\begin{document}

\begin{frontmatter}

\title{Lyapunov-based Stability of Feedback Interconnections of Negative Imaginary Systems} 

\author[First]{Ahmed~G.~Ghallab,}
\author[Second]{Mohamed A. Mabrok,  }
\author[First]{and Ian~R.~Petersen}
\address[First]{School of Engineering and Information Technology, University of New South Wales at the Australian Defence Force
Academy, Canberra ACT 2600, Australia,
email: ahmed.ibrahim@student.adfa.edu.au,  m.a.mabrok@gmail.com, i.r.petersen@gmail.com}
\address[Second]{Robotics, Intelligent Systems and  Control Lab at
Computer, Electrical and Mathematical Science and Engineering,
King Abdullah University of Science and Technology (KAUST)}


\begin{abstract}
Feedback control systems using sensors and actuators such as piezoelectric sensors and actuators, micro-electro-mechanical systems (MEMS) sensors and opto-mechanical sensors, are allowing new advances in designing such high precision technologies. The negative imaginary control systems framework allows for  robust control design for such high precision systems in the face of uncertainties due to un modelled dynamics.  The stability of the feedback interconnection of negative imaginary systems has been well established in the literature. However, the proofs of stability feedback interconnection which are used in some previous papers have a shortcoming due to a matrix inevitability issue.
In this paper we provide a new and correct Lyapunobv-based proof of one such result and show that the result is still true.

\end{abstract}

\end{frontmatter}

\section{Introduction}
Technologies such as atomic force microscopy, nano-positioning, micro-robotics and hard disc drives require high precision and performance in controller design. Feedback control systems using sensors and actuators such as piezoelectric sensors and actuators, micro-electro-mechanical systems (MEMS) sensors and opto-mechanical sensors, \cite{Wilson2002,Harigae20,Bhikkaji2009,Mahmood2011,Salapaka2002,Hulzen2010,Devasia2007,Diaz2012,Ray1978281}, are allowing new advances in designing such high precision technologies. Thew negative imaginary (NI) control systems framework allows for robust control design for such high precision systems in the face of uncertainties due to unmodelled dynamics along with sensor and actuator, \cite{petersen2010,Mabrok2013}.
\begin{figure}
  \centering\includegraphics[width=9cm]{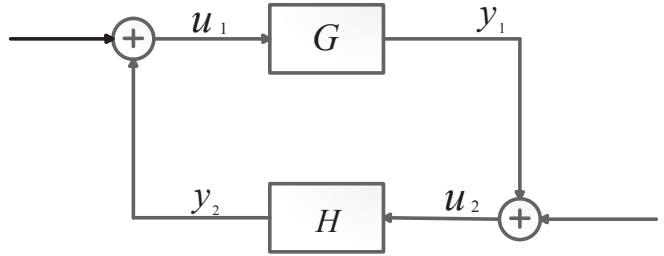}\\
  \caption{A negative-imaginary feedback control system. If the plant
transfer function matrix $G(s)$ is NI and the controller transfer
function matrix $H(s)$ is SNI, then the positive-feedback
interconnection is internally stable if and only if the DC gain
condition, $\lambda_{max}(G(0)H(0))<1,$ is satisfied.
}\label{conn:NI:SNI}
\end{figure}
In  \cite{lanzon2008,petersen2010}, they introduced the theory of negative imaginary (NI) systems for the robust control of flexible structures. Negative imaginary systems are defined by considering the properties of the imaginary part of the system frequency response $G(j\omega)$ and requiring the condition $j(G(j\omega)-G(j\omega)^{\ast})\geq 0$ for all $\omega \in (0,\infty)$.

The robust stability of feedback interconnections of linear time-invariant (LTI) multiple-input multiple-output negative imaginary systems has been studied in \cite{lanzon2008,petersen2010,xiong21010jor}. It is
shown that a necessary and sufficient condition for the internal
stability of a positive-feedback control system  (see Fig.
\ref{conn:NI:SNI}) consisting of an NI plant  with transfer function
matrix $G(s)$ and a strictly negative imaginary (SNI) controller
with transfer function matrix $H(s)$ is given by the DC gain
condition
\begin{equation}\label{DC:ian:alex:cond}
    \lambda_{max}(G(0)H(0))<1,
\end{equation}
 where the notation $\lambda_{max}(\cdot)$ denotes the maximum
eigenvalue of a matrix with only real eigenvalues.  In \cite{Mabrok2013}, a new NI
system definition is given, which allows for flexible structure
systems with colocated force actuators and position sensors, and
with free body motion. Also,   necessary and sufficient
conditions are provided for the stability of  positive feedback
control systems where the plant is NI according to the new
definition and the controller is strictly negative imaginary is given in \cite{Mabrok2013}.

The NI stability result provided in \cite{lanzon2008,petersen2010}
result has been used in a number of  practical applications
\cite{hagen2010,Bhikkaji2009,Mahmood2011,Ahmed2011,Bhikkaji2012,Diaz2012}.
For example in \cite{hagen2010},    this   stability result is
applied to the problem of decentralized control of large vehicle
platoons. In \cite{Bhikkaji2009,Mahmood2011}, the NI stability
result is applied to  nanopositioning  in   an  atomic force
microscope. A positive position feedback control scheme based on the
NI stability result provided in \cite{lanzon2008,petersen2010} is
used to design a novel compensation method for a coupled
fuselage-rotor mode of a rotary wing unmanned aerial vehicle  in
\cite{Ahmed2011}. In \cite{Diaz2012}, an IRC scheme based on the NI stability result  is used
to design  an active vibration control  system for the mitigation of
human induced vibrations in light-weight civil engineering
structures, such as floors and footbridges via proof-mass actuators.
 An identification algorithm which enforces the
NI constraint is proposed in  \cite{Bhikkaji2012} for estimating
model parameters, following which  an Integral resonant controller
is designed for damping vibrations in  flexible structures. In
addition, it is shown in \cite{Schaft2011} that the class of linear
systems having NI transfer function matrices is closely related to
the class of linear port-Hamiltonian input-output systems. Also, an
extension of the NI systems theory  to infinite-dimensional systems
is presented in \cite{Opmeer2011}.

The stability of the feedback interconnection of negative imaginary systems is  established in \cite{lanzon2008,petersen2010,xiong21010jor,Mabrok2013}. However, the proofs of the stability  results which  used in \cite{xiong21010jor, Mabrok2013} have a shortcoming due to a matrix inevitability issue for the case in which the plant has poles on the imaginary axis. In this paper, we are use a Lyapunov-based  stability approach to prove  the result of \cite{xiong21010jor} and show that the result is still true. Note that the result of \cite{KPR1a} does provide a correct proof in the case of plant poles on imaginary axis but requires an extra condition in the definition of the NI property that a certain residue matrix is positive-definite. That extra condition is not required in this paper.

\section{Preliminaries}\label{sec:preliminaries}

In this section, we introduce the concept of negative imaginary systems and strictly negative imaginary systems which allow for poles on the imaginary axis except at the origin.
%
%
%
%
%
\begin{definition}\cite{lanzon2007,lanzon2008,xiong2009}
A square transfer function matrix $G(s)$ is called negative imaginary (NI) if the following conditions are satisfied:
\begin{enumerate}
\item[1.] $G(s )$ has no pole at the origin and in $\Re[s]>0$.
\item[2.] For all $\omega >0$, such that $j\omega$ is not a pole of $G(s )$, $j\left( G(j\omega )-G(j\omega )^{\ast}\right) \geq 0$.
\item[3.] If $j\omega_{0}$, $\omega_0\in(0,\infty)$, is a pole of $G(j\omega )$, it is at most a simple pole and the residue matrix $K_{0}= \lim_{ s\rightarrow j\omega_{0}}(s-j\omega_{0})sG(s)$ is positive semidefinite Hermitian.
\end{enumerate}
\end{definition}

A linear time invariant system is NI if its transfer function is NI.

\begin{definition} \cite{lanzon2008, petersen2010, xiong21010jor}. A square real-rational transfer function matrix $H(s)$ is strictly negative imaginary if:
\begin{enumerate}
\item [1)] $G(s)$ has no poles in $\Re[s]\geq0$;
\item [2)] $j[G(j\omega)-G^{\ast}(j\omega)]>0$ for $\omega\in(0,\infty)$.
\end{enumerate}

\end{definition}

\begin{lemma} \cite{xiong2009}.
\label{PR-NI-lemma}
Given a real rational strictly proper transfer function matrix
$G(s)$ with minimal state space realization $
\begin{bmatrix}
\begin{array}{c|c}
A & B \\ \hline C & D
\end{array}
\end{bmatrix}
$ and, define  the transfer function matrix $\tilde{G}(s)=G(s)-D $.
The transfer function matrix ${G}(s)$ is negative imaginary if and
only if,

\begin{enumerate}
\item ${G}(s)$ has no poles at the origin.

\item The transfer function matrix $F(s)=s\tilde{G}(s)$ is positive real.
\end{enumerate}
\end{lemma}

\begin{lemma}\cite{lanzon2008}. \label{NI}
Let $(A,B,C,D)$ be a minimal state-space realization of the $m\times m$ real-rational proper transfer function matrix $G(s)$, where $A \in \mathbb{R}^{n \times n},B \in \mathbb{R}^{n \times m},C\in \mathbb{R}^{m \times n},D \in \mathbb{R}^{m \times m}$. Then $G(s)$ is negative imaginary if and only if:
\begin{enumerate}
\item [1)] det$(A)\neq 0$, $D=D^{T}$;
\item [2)] there exists a matrix $P=P^{T}>0, P \in \mathbb{R}^{n \times n}$, such that
\end{enumerate}
$$
AP^{-1}+P^{-1}A^{T}\leq 0, \quad \text{and } \quad B+AP^{-1}C^{T}=0
$$
\end{lemma}

\begin{lemma}\label{SNI}\cite{xiong21010jor}
Let $(A,B,C,D)$ be a minimal state-space realization of the $m\times m$ real-rational proper transfer function matrix $G(s)$, where $A \in \mathbb{R}^{n \times n},B \in \mathbb{R}^{n \times m},C\in \mathbb{R}^{m \times n},D \in \mathbb{R}^{m \times m}$. Then $G(s)$ is strictly negative imaginary if and only if:
\begin{enumerate}
\item [1)] det$(A)\neq 0$, $D=D^{T}$;
\item [2)] there exists a matrix $P=P^{T}>0, P \in \mathbb{R}^{n \times n}$, such that
$$
AP^{-1}+P^{-1}A^{T}\leq 0, \quad \text{and } \quad B+AP^{-1}C^{T}=0
$$
\item [3)]the transfer function matrix $M(s)\backsim
\begin{bmatrix}
\begin{array}{c|c}
A & B \\ \hline LPA^{-1} & 0
\end{array}
\end{bmatrix}$
 has full column rank at $s=jw$ for any $\omega\in(0,\infty)$ where $L^{T}L=-AP^{-1}-P^{-1}A^{T}$. That is, rank $M(j\omega)=m$ for any $\omega\in(0,\infty)$.
\end{enumerate}
\end{lemma}

We will consider the positive feedback interconnection of a linear NI system with a linear SNI system and prove the internal stability of the closed-loop system. Consider a minimal state-space representation for the SNI transfer function $H(s)$,
\begin{align}
\label{eq:xdot1}
&\dot{x_{1}}(t) = A_{1} x_{1}(t)+B_{1} u_{1}(t), \\
\label{eq:y1} &y_{1}(t) = C_{1} x_{1}(t)+D_{1} u_{1}(t),
\end{align}%
where $A_2 \in \mathbb{R}^{n \times n},B_2 \in \mathbb{R}^{n \times m},C_2
\in \mathbb{R}^{m \times n},D_2 \in \mathbb{R}^{m \times m}$.

Also, we consider a minimal state-space representation for the NI transfer function $G(s)$,
\begin{align}
\label{eq:xdot2}
&\dot{x_{2}}(t) = A_{2} x_{2}(t)+B_{2} u_{2}(t), \\
\label{eq:y2} &y_{2}(t) = C_{2} x_{2}(t)+D_{2} u_{2}(t),
\end{align}%
where $A_1 \in \mathbb{R}^{n \times n},B_1 \in \mathbb{R}^{n \times m},C_1
\in \mathbb{R}^{m \times n},D_1 \in \mathbb{R}^{m \times m}$.

Since $G(s)$ is NI, Lemma \ref{NI} implies there exist a symmetric matrix $P_{1}>0$, and a matrix $L_{1}$ such that

\begin{align*}
  A_{1}P_{1}^{-1}+A_{1}^{T}P_{1}^{-1}&= L^{T}_{1}L_{1} \\
 B_{1}+A_{1}P_{1}^{-1}C_{1}^{T}&=0 \\
\end{align*}
which give
\begin{align*}
 P_{1}A_{1}+A_{1}^{T}P_{1}&= -P_{1}L^{T}_{1}L_{1}P_{1} \\
 B_{1}^{T}P_{1}-C_{1}L^{T}_{1}L_{1}P_{1}&= C_{1}A_{1}\\
 C_{1}B_{1}+(C_{1}B_{1})^{T}&=(L_{1}C_{1}^{T})^{T}(L_{1}C_{1}^{T})
\end{align*}

Also, since $H(s)$ is SNI, Lemma \ref{SNI} implies there exist a symmetric matrix $P_{2}>0$, and a matrix $L_{2}$ such that

\begin{align*}
 A_{2}P_{2}^{-1}+A_{2}^{T}P_{2}^{-1}&= L_{2}^{T}L_{2} \\
 B_{2}+A_{2}P_{2}^{-1}C_{2}^{T}&=0 \\
\end{align*}
which gives
\begin{align*}
P_{2}A_{2}+A_{2}^{T}P_{2}&= -P_{2}L^{T}_{2}L_{2}P_{2} \\
B_{2}^{T}P_{2}-C_{2}L^{T}_{2}L_{2}P_{2}&= C_{2}A_{2}\\
C_{2}B_{2}+(C_{2}B_{2})^{T}&=(L_{2}C_{2}^{T})^{T}(L_{2}C_{2}^{T})
\end{align*}

where $$W(s):=sM(s)=L_2P_2(sI-A_2)^{-1}B_2-L_2C_2^{T}$$ has no zeros in the $j\omega$-axis except at the origin.

\begin{lemma}\label{lemma1}
Given negative imaginary $G(s)$ and strictly negative imaginary $H(s)$. Assume $G(\infty)H(\infty)=0$ and $N(\infty)\geq0$. Then the matrix
$$
\left[ \begin{array}{cc}
{P}_{1}-C_{1}^{T}D_{2}C_{1} & -C_{1}^{T}C_{2} \\
{-C}_{2}^{T}C_{1} & {P}_{2}-C_{2}^{T}D_{1}C_{2} \end{array}
\right]
$$
is positive definite if and only if $\lambda_{max}(G(0)H(0))<1$.

\begin{pf} See also \cite{lanzon2008}. We have
\begin{align*}
     \lambda&_{max}(G(0)H(0))<1\\[8pt]
     & \Leftrightarrow H(0)^{-1}-G(0)>0\\[8pt]
     & \Leftrightarrow H(0)^{-1}-D-C_{1}P_{1}C_{1}^{T}>0\\[8pt]
     & \Leftrightarrow
     \left[ \begin{array}{cc}
      {P}_{1} & -C_{1}^{T}\\
      C_{1} & H(0)^{-1}-D_{1} \end{array}
     \right]\\[8pt]
     & \Leftrightarrow H(0)^{-1}-D>0, \quad \text{and }\\[8pt]
     & \quad P_{1}-C_{1}(H(0)^{-1}-D_{1})^{-1}C_{1}>0\\[8pt]
     & \Leftrightarrow \lambda_{max}[D_1H(0)]<1, \quad \text{and } \\[8pt]
     &\quad P_{1}-C_{1}(H(0)^{-1}-D_{1})^{-1}[D_2+(H(0)-D_2)]C_{1}>0\\[8pt]
     & \Leftrightarrow\lambda_{max}[D_{1}C_{2}P_{2}^{-1}C_{2}^{T}], \quad \text{and }\\[8pt]
     & \quad P_{1}-C_{1}^{T}D_{2}C_{1}-C_{1}^{T}(I-H(0)D_{1})^{-1}(H(0)-D_{2})C_{1}>0\\[8pt]
     &\Leftrightarrow\lambda_{max}[P_{2}^{\frac{1}{2}}C_{2}^{T}D_{1}C_{2}P_{2}^{\frac{1}{2}}], \quad \text{and }\\[8pt]
     &\quad   P_{1}-C_{1}^{T}D_{2}C_{1}-C^{T}(I-C_{2}P_{2}C_{2}^{T}D_{1})C_{2}P_{2}C_{2}^{T}C_{1}>0\\[8pt]
     &\Leftrightarrow P_{2}-C_{2}^{T}D_{1}C_{2}>0, \quad \text{and }\\[8pt]
     &\quad  (P_{1}-C_{1}^{T}D_{2}C_{1})-C_{1}^{T}C_{2}(P_{2}-C_{2}^{T}D_{1}C_{2})^{-1}C_{2}^{T}C_{1}>0\\[8pt]
     &\Leftrightarrow
     \left[ \begin{array}{cc}
          {P}_{1}-C_{1}^{T}D_{2}C_{1} & -C_{1}^{T}C_{2} \\
          {-C}_{2}^{T}C_{1} & {P}_{2}-C_{2}^{T}D_{1}C_{2} \end{array}
            \right]>0.
   \end{align*}
\begin{flushright}
$\square$
\end{flushright}
\end{pf}
\end{lemma}

\section{Main results}\label{sec:main-results}

In this section we shall introduce the main result regarding the internal stability of the positive feedback interconnection of $G(s)$ and $H(s)$.
\begin{theorem}
\label{th:ni-thm} Assume $G(s)$ is negative imaginary system and $H(s)$ is strictly negative imaginary system such that  $G(\infty)H(\infty)=0$ and $N(\infty)\geq0$. Also, assume that $\lambda_{max}(G(0)H(0))<1$. Then, the positive feedback interconnection of $G(s)$ and $H(s)$ is internally stable.
\end{theorem}

\begin{pf} This proof follows a similar approach to the proof of Lemma 3.37 in \cite{brogliato-bk2007}.

\noindent
Let $V_1(x_1)=x_{1}^{T}P_1x_1$ and $V_2(x_2)=x_{2}^{T}P_2x_2$ and consider the function
\begin{align*}
V(x_{1},x_{2})&= V_1(x_1)+V_2(x_2)-2y_{1}^{T}y_{2}\\[8pt]
              &=  \left[\begin{array}{r}
                x_{1}^{T} \hspace{1.5mm} x_{2}^{T}
               \end{array}\right]
               \left[ \begin{array}{cc}
               {P}_{1}-C_{1}^{T}D_{2}C_{1} & -C_{1}^{T}C_{2} \\
               {-C}_{2}^{T}C_{1} & {P}_{2}-C_{2}^{T}D_{1}C_{2} \end{array}
                \right]
               \left[ \begin{array}{c}
                  x_{1} \\
                  x_{2} \\
               \end{array}\right]
\end{align*}
as a Lyapunov candidate for the feedback system. Note that it follows from Lemma \ref{lemma1} that the function $V(x_1,x_2)$ is positive definite. Now for the closed loop system we have
\begin{align*}
&\dot{V}(x_{1},x_{2})=\\
&\left[\begin{array}{r}
\dot{x}_{1}^{T} \hspace{1.5mm} \dot{x}_{2}^{T}
\end{array}\right]
\left[ \begin{array}{cc}
{P}_{1}-C_{1}^{T}D_{2}C_{1} & -C_{1}^{T}C_{2} \\
{-C}_{2}^{T}C_{1} & {P}_{2}-C_{2}^{T}D_{1}C_{2} \end{array}
\right]
\left[ \begin{array}{c}
x_{1} \\
x_{2} \\
\end{array}\right]\\[8pt]
&\quad + \left[\begin{array}{r}
\dot{x}_{1}^{T} \hspace{1.5mm} \dot{x}_{2}^{T}
\end{array}\right]
\left[ \begin{array}{cc}
{P}_{1}-C_{1}^{T}D_{2}C_{1} & -C_{1}^{T}C_{2} \\
{-C}_{2}^{T}C_{1} & {P}_{2}-C_{2}^{T}D_{1}C_{2} \end{array}
\right]
\left[ \begin{array}{c}
\dot{x}_{1} \\
\dot{x}_{2} \\
\end{array}\right]\\[8pt]
&=x_{1}^{T}P_{1}\dot{x}_{1}+\dot{x}_{1}^{T}P_{1}x_{1}-2\dot{x}_{1}^{T}C_{1}^{T}D_{2}C_{1}x_{1}-2x_{2}^{T}C_{2}^{T}C_{1}x_{1}\\[8pt]
&\quad -2x_{1}^{T}C_{1}^{T}C_{2}x_{2}+\dot{x}_{2}^{T}P_{2}x_{2}+x_{2}^{T}P_{2}\dot{x}_{2}-2x_{2}^{T}C_{2}^{T}D_{1}C_{2}\dot{x}_{2}\\[8pt]
&=(x_{1}^{T}A_{1}^{T}+u_{1}^{T}B_{1}^{T})P_1x_1+x_{1}^{T}P_1(A_1x_1+B_1u_1)\\[8pt]
&\quad + (x_{2}^{T}A_{2}^{T}+u_{2}^{T}B_{2}^{T})P_2x_2+x_{2}^{T}P_2(A_2x_2+B_2u_2)\\[8pt]
&\quad  -2(\dot{y}_{1}^{T}-\dot{u}_{1}^{T}D_{1}^{T})D_2(y_1-D_1u_1)\\[8pt]
&\quad -2(\dot{y}_{2}^{T}-\dot{u}_{2}^{T}D_{2}^{T})(y_1-D_1u_1)\\[8pt]
&\quad - 2(\dot{y}_{1}^{T}-\dot{u}_{1}^{T}D_{1}^{T})(y_2-D_2u_2)\\[8pt]
&\quad -2(\dot{y}_{2}^{T}-\dot{u}_{2}^{T}D_{2}^{T})D_1(y_2-D_2u_2)\\[8pt]
&=x_{1}^{T}(A_{1}^{T}P_{1}+P_{1}A_{1})x_1+x_{2}^{T}(A_{2}^{T}P_{2}+P_{2}A_{2})x_2\\[8pt]
&\quad+2u_{1}^{T}B_{1}^{T}P_{1}x_{1}+2u_{2}^{T}B_{2}^{T}P_{2}x_{2}\\[8pt]
&\quad-2\dot{y}_{1}^{T}D_2y_1-2\dot{y}_{2}^{T}y_1+2\dot{y}_{2}^{T}D_1u_1+2\dot{u}_{2}^{T}D_{2}^{T}y_1\\[8pt]
&\quad+2\dot{y}_{1}^{T}D_2u_2+2\dot{u}_{1}^{T}D_{1}^{T}y_{2}-2\dot{y}_{2}^{T}D_1y_2-2\dot{y}_{1}^{T}y_2\\[8pt]
&=-x_{1}^{T}P_{1}L^{T}_{1}L_{1}P_{1}x_{1}-x_{2}^{T}P_{2}L^{T}_{2}L_{2}P_{2}x_{2}-2\dot{y}_{2}^{T}y_{1}-2\dot{y}_{1}^{T}y_2\\[8pt]
&\quad+2u^{T}_{1}(B_{1}^{T}P_{1}-C_{1}L_{1}^{T}L_{1}P_{1})x_{1}-2u^{T}_{1}C_{1}L_{1}^{T}L_{1}P_{1}x_{1}\\[8pt]
&\quad +2u^{T}_{2}(B_{2}^{T}P_{2}-C_{2}L_{2}^{T}L_{2}P_{2})x_{2}-2u^{T}_{2}C_{2}L_{2}^{T}L_{2}P_{2}x_{2}\\[8pt]
&= -x_{1}^{T}P_{1}L^{T}_{1}L_{1}P_{1}x_{1}-x_{2}^{T}P_{2}L^{T}_{2}L_{2}P_{2}x_{2}-2\dot{y}_{2}^{T}y_{1}-2\dot{y}_{1}^{T}y_2\\[8pt]
&\quad +2u^{T}_{1}(C_{1}A_{1}x_{1}+C_{1}B_{1}u_{1})-2u^{T}_{1}C_{1}B_{1}u_{1}\\[8pt]
&\quad +2u^{T}_{2}(C_{2}A_{2}x_{2}+C_{2}B_{2}u_{2})-2u^{T}_{2}C_{2}B_{2}u_{2}\\[8pt]
&\quad -2u^{T}_{1}C_{1}L_{1}^{T}L_{1}P_{1}x_{1}-2u^{T}_{2}C_{2}L_{2}^{T}L_{2}P_{2}x_{2}\\[8pt]
&= -x_{1}^{T}P_{1}L^{T}_{1}L_{1}P_{1}x_{1}-x_{2}^{T}P_{2}L^{T}_{2}L_{2}P_{2}x_{2}-2\dot{y}_{2}^{T}y_{1}-2\dot{y}_{1}^{T}y_2\\[8pt]
&\quad +2u^{T}_{1}\dot{y}_{1}-2u^{T}_{1}C_{1}B_{1}u_{1}-2u^{T}_{1}C_{1}L_{1}^{T}L_{1}P_{1}x_{1}\\[8pt]
&\quad +2u^{T}_{2}\dot{y}_{2}-2u^{T}_{2}C_{2}B_{2}u_{2}-2u^{T}_{2}C_{2}L_{2}^{T}L_{2}P_{2}x_{2}\\[8pt]
&= -(L_{1}P_{1}x_{1}-L_{1}C^{T}_{1}u_{1})^{T}(L_{1}P_{1}x_{1}-L_{1}C^{T}_{1}u_{1})\\[8pt]
&\quad \quad \quad \quad-(L_{2}P_{2}x_{2}-L_{2}C^{T}_{2}u_{2})^{T}(L_{2}P_{2}x_{2}-L_{2}C^{T}_{2}u_{2})
\end{align*}

where we used  the equations $u_{1}=y_{2}$ and $u_{2}=y_{1}$, and 
\begin{equation*}
2u^{T}_{i}C_{i}B_{i}u_{i}= u^{T}_{i}(C_{i}B_{i}+(C_{i}B_{i})^{T})u_{i}=u^{T}_{i}C_{i}L_{i}^{T}L_{i}C_{i}^{T}u_{i}.
\end{equation*}

 Define $\tilde{y}_{i}=L_{i}P_{i}x_{i}-L_{i}C^{T}_{i}u_{i}$, for $i=1,2$. Then
\begin{equation}\label{eq:2}
\dot{V}(x_1,x_2)= -\tilde{y}_{1}^{T}\tilde{y}_{1} -\tilde{y}_{2}^{T}\tilde{y}_{2}\leq -\tilde{y}_{2}^{T}\tilde{y}_{2}\leq 0.
\end{equation}

This implies that the closed loop systems is at least Lyapunov-stable, \ie \ the closed loop system poles are only in the closed left half plane. 

Integrating \eqref{eq:2} we get
\begin{equation}\label{eq:1}
-V(0)\leq V(t)-V(0)\leq - \int_{0}^{t} \tilde{y}_{2}^{T}(s)\tilde{y}_{2}(s)ds,
\end{equation}
then
\begin{equation}\label{eq:3}
  \int_{0}^{t}\tilde{y}_{2}^{T}(s)\tilde{y}_{2}(s)ds \leq V(0),
\end{equation}
which implies
\begin{equation}\label{eq:22}
\tilde{y}_{2}=L_{2}P_{2}x_{2}-L_{2}C^{T}_{2}u_{2}=0.
\end{equation}

We now show that the closed loop system matrix has no eigen values on the imaginary axis. The closed loop matrix for the systems is
$$\breve{A}=
 \left[ \begin{array}{cc}
          {A}_{1}+B_{1}D_{2}C_{1} & B_{1}C_{2} \\
          {B}_{2}C_{1} & {A}_{2}+B_{2}D_{1}C_{2} \end{array}
\right].
$$

Suppose that this matrix has an eigen value on the $j\omega$-axis. Then there exists a nonzero $x=\left[\begin{array}{r}
           x_{1}^{T} \hspace{1.5mm} x_{2}^{T}
          \end{array}\right]^{T}$ such that
$$
 \left[ \begin{array}{cc}
          {A}_{1}-j\omega I+B_{1}D_{2}C_{1} & B_{1}C_{2} \\
          {B}_{2}C_{1} & {A}_{2}-j\omega I+B_{2}D_{1}C_{2} \end{array}
\right]
\left[\begin{array}{c}
                  x_{1} \\
                  x_{2} \\
               \end{array}\right]=0,
$$

for $\omega\in\mathbb{R}$. So, we have
\begin{equation}\label{ee1}
 (A_1-j\omega I+B_1D_2C_1)x_1+B_1C_2x_2 = 0,
\end{equation}
and
\begin{equation}\label{ee2}
B_2C_1x_1+(A_2-j\omega I+B_2D_1C_2)x_2= 0.
\end{equation}
Then, we have
\begin{equation}\label{eq:12}
 (j\omega I-A_{1})x_{1}-B_{1}y_{2}=0,
\end{equation}
and
\begin{equation}\label{eq:21}
  (j\omega I-A_{2})x_{2}-B_{2}y_{1}=0,
\end{equation}
where we used the state-space equations \eqref{eq:y1}, \eqref{eq:y2} and the equations $u_{1}=y_{2}$ and $u_{2}=y_{1}$.

Combining  \eqref{eq:21}, \eqref{eq:22} in a matrix equation form we get
$$
 \left[ \begin{array}{cc}
          {A}_{2}-j\omega I & B_{2} \\
          {L}_{2}P_{2} & -{L}_{2}C_{2}^{T} \end{array}
\right]\left[ \begin{array}{c}
                  x_{2} \\
                  u_{2} \\
               \end{array}\right]=0.
$$

Since the matrix on the left has full rank for $\omega\in(0,\infty)$, it follows that $x_{2}=u_{2}=0$ and hence $y_1=y_2=0$. This implies $C_{1}x_{1}=0, x_{1}\neq0$, \ie \ $(A,C)$ is non-observable, which contradicts the minimality of the systems $G(s)$. Therefore, $\breve{A}$ is semistable (\ie $\ j\omega \notin \text{spec}(\breve{A}), \ \text{for nonzero} \ \omega\in \mathbb{R}$).

Now assume that the matrix $\breve{A}$ has an eigen value at the origin ($\omega=0$). From \eqref{eq:12} we have
$$
C_1x_1=-C_1A_{1}^{-1}B_{1}y_2=(G(0)-D_1)y_2.
$$
This implies
$$
y_1-D_1u_1=(G(0)-D_1)y_2,
$$
and then
\begin{equation}\label{eq:102}
y_1=G(0)y_{2} .
\end{equation}
Similarly, from \eqref{eq:21} we have
\begin{equation}\label{eq:201}
y_2=H(0)y_{1}.
\end{equation}
Combining \eqref{eq:102} and \eqref{eq:201} we get
\begin{equation}\label{1001}
  y_1=G(0)H(0)y_1.
\end{equation}

Note that from \eqref{eq:21}, if $y_1=0$ we get $x_2=0$, since $A_2$ is asymptotically stable and hence invertible. Also, we have $y_2=C_2x_2+D_2u_2=0,$ and from \eqref{eq:12} we have $x_1=0$ since $A_1$ has no have eigenvalue at the origin. That leads to $(x_1,x_2)=0$ which is not allowed, thus $y_1$ must be nonzero. However, \eqref{1001} contradicts with the DC gain condition $\lambda_{max}(G(0)H(0))<1$. Therefore, we have shown by contradiction that the closed loop system does not have eigen values on the imaginary axis. From that we conclude that the feedback interconnection of $G(s)$ and $H(s)$ is internally stable.
\begin{flushright}
$\square$
\end{flushright}



\end{pf}
The next corollary shows that the result in \cite{xiong21010jor} for the internal stability of positive feedback interconnections of NI systems is still true.
\begin{cor}
Given a NI transfer function matrix $G(s)$ and a SNI transfer function matrix $H(s)$ and assume that $G(\infty)H(\infty)=0$ and $H(\infty)\geq 0$. Then, the feedback interconnection of $G(s)$ and $H(s)$ is internally stable if and only if $\lambda_{max}(G(0)H(0))<1$.
\end{cor}
\begin{pf}
This result follows from Theorem 1 and the necessity part of Theorem 5 of \cite{xiong21010jor} for which a correct proof has already been given in \cite{xiong21010jor}.
\begin{flushright}
$\square$
\end{flushright}
\end{pf}

\section*{Acknowledgment}
The authors would like to acknowledge Sei Zhen Khong who originally pointed out to the problem with the proof of Theorem 5 in \cite{xiong21010jor}. This work was supported by the Australian Research Council under grant DP160101121.

\balance



\begin{thebibliography}{22}
\providecommand{\natexlab}[1]{#1}
\providecommand{\url}[1]{\texttt{#1}}
\expandafter\ifx\csname urlstyle\endcsname\relax
  \providecommand{\doi}[1]{doi: #1}\else
  \providecommand{\doi}{doi: \begingroup \urlstyle{rm}\Url}\fi

\bibitem[Ahmed and Pota(2011)]{Ahmed2011}
B.~Ahmed and H.~Pota.
\newblock Dynamic compensation for control of a rotary wing {UAV} using
  positive position feedback.
\newblock \emph{Journal of Intelligent and Robotic Systems}, 61\penalty0
  (1-4):\penalty0 43--56, 2011.
\newblock ISSN 0921-0296.

\bibitem[Bhikkaji and Moheimani(2009)]{Bhikkaji2009}
B.~Bhikkaji and S.~Moheimani.
\newblock Fast scanning using piezoelectric tube nanopositioners: A negative
  imaginary approach.
\newblock In \emph{Proc. IEEE/ASME Int. Conf. Advanced Intelligent Mechatronics
  AIM}, pages 274--279, Singapore, July 2009.
\newblock \doi{10.1109/AIM.2009.5230001}.

\bibitem[Bhikkaji et~al.(2012)Bhikkaji, Moheimani, and Petersen]{Bhikkaji2012}
B.~Bhikkaji, S.~O.~R. Moheimani, and I.~R. Petersen.
\newblock A negative imaginary approach to modeling and control of a collocated
  structure.
\newblock \emph{IEEE/ASME Transactions on Mechatronics}, 17\penalty0
  (4):\penalty0 717--727, 2012.
\newblock \doi{10.1109/TMECH.2011.2123909}.

\bibitem[Brogliato et~al.(2007)Brogliato, Lozano, Maschke, and
  Egeland]{brogliato-bk2007}
B.~Brogliato, R.~Lozano, B.~Maschke, and O.~Egeland.
\newblock \emph{Dissipative Systems Analysis and Control}.
\newblock Communications and Control Engineering. Springer, London, UK, 2nd
  edition, 2007.

\bibitem[Cai and Hagen(Aug. 2010)]{hagen2010}
C.~Cai and G.~Hagen.
\newblock Stability analysis for a string of coupled stable subsystems with
  negative imaginary frequency response.
\newblock \emph{IEEE Transactions on Automatic Control}, 55\penalty0
  (8):\penalty0 1958--1963, Aug. 2010.

\bibitem[D.~G.~Wilson and Starr(2002)]{Wilson2002}
G.~G.~Parker D.~G.~Wilson, R. D.~Robinett and G.~P. Starr.
\newblock Augmented sliding mode control for flexible link manipulators.
\newblock \emph{Journal of Intelligent and Robotic Systems}, 34\penalty0
  (4):\penalty0 415--430, 2002.

\bibitem[Devasia et~al.(2007)Devasia, Eleftheriou, and Moheimani]{Devasia2007}
S.~Devasia, E.~Eleftheriou, and S.~O.~R. Moheimani.
\newblock A survey of control issues in nanopositioning.
\newblock \emph{IEEE Transactions on Control Systems Technology}, 15\penalty0
  (5):\penalty0 802--823, 2007.
\newblock \doi{10.1109/TCST.2007.903345}.

\bibitem[Diaz et~al.(2012)Diaz, Pereira, and Reynolds]{Diaz2012}
I.~M. Diaz, E.~Pereira, and P.~Reynolds.
\newblock Integral resonant control scheme for cancelling human-induced
  vibrations in light-weight pedestrian structures.
\newblock \emph{Structural Control and Health Monitoring}, 19\penalty0
  (1):\penalty0 55--69, 2012.

\bibitem[Harigae et~al.(2003)Harigae, Yamaguchi, Kasai, Igawa, and
  Suzuki]{Harigae20}
M.~Harigae, I.~Yamaguchi, T.~Kasai, H.~Igawa, and T.~Suzuki.
\newblock Control of large space structures using {GPS} modal parameter
  identification and attitude and deformation estimation.
\newblock \emph{Electronics and Communications in Japan}, 86\penalty0
  (4):\penalty0 63--71, 2003.

\bibitem[Lanzon and Petersen(2007)]{lanzon2007}
A.~Lanzon and I.~R. Petersen.
\newblock A modified positive-real type stability condition.
\newblock In \emph{Proceedings of the European Control Conference}, pages
  3912--3918, Kos, Greece, July 2007.

\bibitem[Lanzon and Petersen(2008)]{lanzon2008}
A.~Lanzon and I.~R. Petersen.
\newblock Stability robustness of a feedback interconnection of systems with
  negative imaginary frequency response.
\newblock \emph{IEEE Transactions on Automatic Control}, 53\penalty0
  (4):\penalty0 1042--1046, 2008.

\bibitem[M.~A.~Mabrok and Lanzon(2014)]{Mabrok2013}
I.~R.~Petersen M.~A.~Mabrok, A. G.~Kallapur and A.~Lanzon.
\newblock Generalizing negative imaginary systems theory to include free body
  dynamics: Control of highly resonant structures with free body motion.
\newblock \emph{Automatic Control, IEEE Transactions on}, 59\penalty0
  (10):\penalty0 2692--2707, Oct 2014.
\newblock ISSN 0018-9286.
\newblock \doi{10.1109/TAC.2014.2325692}.

\bibitem[Mahmood et~al.(2011)Mahmood, Moheimani, and Bhikkaji]{Mahmood2011}
I.~A. Mahmood, S.~O.~R. Moheimani, and B.~Bhikkaji.
\newblock A new scanning method for fast atomic force microscopy.
\newblock \emph{IEEE Transactions on Nanotechnology}, 10\penalty0 (2):\penalty0
  203--216, 2011.

\bibitem[Opmeer(2011)]{Opmeer2011}
M.~R. Opmeer.
\newblock Infinite-dimensional negative imaginary systems.
\newblock \emph{IEEE Transactions on Automatic Control}, 56\penalty0
  (12):\penalty0 2973--2976, 2011.
\newblock \doi{10.1109/TAC.2011.2162886}.

\bibitem[Petersen and Lanzon(2010)]{petersen2010}
I.~R. Petersen and A.~Lanzon.
\newblock Feedback control of negative imaginary systems.
\newblock \emph{IEEE Control System Magazine}, 30\penalty0 (5):\penalty0
  54--72, 2010.

\bibitem[Ray(1978)]{Ray1978281}
W.~H. Ray.
\newblock Some recent applications of distributed parameter systems theory-a
  survey.
\newblock \emph{Automatica}, 14\penalty0 (3):\penalty0 281 -- 287, 1978.
\newblock ISSN 0005-1098.

\bibitem[S.~Z.~Khong and Rantzer(2015)]{KPR1a}
I.~R.~Petersen S.~Z.~Khong and A.~Rantzer.
\newblock Robust feedback stability of negative imaginary systems.
\newblock In \emph{Proceedings of the European Control Conference 2015}, Linz,
  Austria, July 2015.

\bibitem[Salapaka et~al.(2002)Salapaka, Sebastian, Cleveland, and
  Salapaka]{Salapaka2002}
S.~Salapaka, A.~Sebastian, J.~P. Cleveland, and M.~V. Salapaka.
\newblock High bandwidth nano-positioner: A robust control approach.
\newblock \emph{Review of Scientific Instruments}, 73\penalty0 (9):\penalty0
  3232--3241, 2002.

\bibitem[van~der Schaft(2011)]{Schaft2011}
A.~J. van~der Schaft.
\newblock Positive feedback interconnection of {H}amiltonian systems.
\newblock In \emph{Proceedings of the 50th IEEE Conference on Decision and
  Control and European Control Conference (CDC-ECC)}, Orlando, FL, USA, Dec
  2011.

\bibitem[van Hulzen et~al.(2010)van Hulzen, Schitter, Van~den Hof, and van
  Eijk]{Hulzen2010}
J.~R. van Hulzen, G.~Schitter, P.~M.~J. Van~den Hof, and J.~van Eijk.
\newblock Modal actuation for high bandwidth nano-positioning.
\newblock In \emph{Proc. American Control Conference}, pages 6525--6530,
  Baltimore, Maryland, USA, July 2010.

\bibitem[Xiong et~al.(2009)Xiong, Petersen, and Lanzon]{xiong2009}
J.~Xiong, I.~R. Petersen, and A.~Lanzon.
\newblock Stability analysis of positive feedback interconnections of linear
  negative imaginary systems.
\newblock In \emph{Proceedings of the American Control Conference}, pages
  1855--1860, St. Louis, Missouri, USA, June 2009.

\bibitem[Xiong et~al.(2010)Xiong, Petersen, and Lanzon]{xiong21010jor}
J.~Xiong, I.~R. Petersen, and A.~Lanzon.
\newblock A negative imaginary lemma and the stability of interconnections of
  linear negative imaginary systems.
\newblock \emph{IEEE Transactions on Automatic Control}, 55\penalty0
  (10):\penalty0 2342--2347, 2010.

\end{thebibliography}
\end{document}